\documentclass[12pt]{article}
\usepackage{amsmath,amsfonts,amsthm,amssymb}
\newcommand{\F}{\mathbb F_q}
\newcommand{\K}{\overline{K}_c}
\newcommand{\AG}{\mathfrak A}
\newcommand{\MG}{\mathfrak M}
\newcommand{\FF}{\mathcal F}
\DeclareMathOperator{\card}{card}

\numberwithin{equation}{section}
\begin{document}
\newtheorem{lem}{Lemma}
\newtheorem{teo}{Theorem}
\pagestyle{plain}
\title{Evolution Equations and Functions of Hypergeometric Type
over Fields of Positive Characteristic}
\author{Anatoly N. Kochubei\footnote{Partially supported by
DFG under Grant 436 UKR 113/87/01, and by the Ukrainian Foundation 
for Fundamental Research, Grant 10.01/004.}\\ 
\footnotesize Institute of Mathematics,\\ 
\footnotesize National Academy of Sciences of Ukraine,\\ 
\footnotesize Tereshchenkivska 3, Kiev, 01601 Ukraine
\\ \footnotesize E-mail: \ kochubei@i.com.ua}
\date{}
\maketitle

\bigskip
\begin{abstract}
We consider a class of partial differential equations with Carlitz 
derivatives over a local field of positive characteristic, for
which an analog of the Cauchy problem is well-posed. Equations of 
such type correspond to quasi-holonomic modules over the  
ring of differential operators with Carlitz derivatives. The above 
class of equations includes some equations of hypergeometric type. 
Building on the work of Thakur, we develop his notion of the 
hypergeometric function of the first kind (whose parameters 
belonged initially to $\mathbb Z$) in such a way that it becomes 
fully an object of the function field arithmetic, with the 
variable, parameters and values from the field of positive 
characteristic.
\end{abstract}



\bigskip
\section{INTRODUCTION}
Let $K$ be the field of formal Laurent series
$t=\sum\limits_{j=N}^\infty \xi_jx^j$ with coefficients $\xi_j$ 
from the Galois field $\F$, $\xi_N\ne 0$ if $t\ne 0$, $q=p^\upsilon $, 
$\upsilon \in \mathbf Z_+$,
where $p$ is a prime number. It is well known that any non-discrete 
locally compact field of characteristic $p$ is isomorphic to such $K$.
The absolute value on $K$ is given by $|t|=q^{-N}$, $|0|=0$. This 
absolute value can be extended in a unique way onto the completion 
$\K$ of an algebraic closure of $K$. See \cite{Sch} for a detailed 
description of the extension procedure (valid for any 
non-Archimedean valued field).

An important feature of analysis over $K$ and $\K$ initiated  
by Carlitz \cite{C35} and developed by Wagner, Goss, 
Thakur, the author, and many others (see surveys in \cite{G,Th3,Ksurv})
is the availability of many non-trivial 
$\F$-linear functions, that is such functions $f$ defined on 
$\F$-subspaces $K_0\subset K$, with values from $\K$, that
$$
f(t_1+t_2)=f(t_1)+f(t_2),\quad f(\alpha t)=\alpha f(t),
$$
for any $t,t_1,t_2\in K_0$, $\alpha \in \F$. Such are, for 
example, polynomials and power series of the form $\sum 
a_kt^{q^k}$. Within this class, there are analogs of the 
exponential and logarithm, the Bessel and hypergeometric 
functions, the polylogarithms and zeta function.

The basic ingredients of the calculus of $\F$-linear functions 
\cite{K99} are the Frobenius operator $\tau u=u^q$, the difference 
operator 
\begin{equation}
\Delta u(t)=u(xt)-xu(t)
\end{equation}
introduced in \cite{C35}, and the nonlinear ($\F$-linear) operator 
$d=\tau^{-1}\Delta$ called {\it the Carlitz derivative}. The latter appears, 
as an analog of the classical derivative, in the theory of ordinary
differential equations for $\F$-linear functions, which has been developed 
both in the traditional analytic direction (the Cauchy problem \cite{K00,Knlin},
regular singularity \cite{K03}, new special functions defined via 
differential equations \cite{Kzeta}) and within various algebraic frameworks
(analogs of the canonical commutation relations of quantum mechanics 
\cite{K98,K99}, umbral calculus \cite{K05}, an analog of the Weyl 
algebra \cite{K00,Khol,Bav}). Note that in our situation the 
meaning of a polynomial coefficient in a differential equation is 
not a usual multiplication by a polynomial, but the action of a 
polynomial in the operator $\tau$ (similarly, for holomorphic 
coefficients). Thus, the operator $\Delta =\tau d$ can be seen as 
the counterpart of $t\dfrac{d}{dt}$.

As an example, consider Thakur's hypergeometric function 
\cite{Th1,Th2,Th3} (in this paper we deal only with what Thakur 
calls ``the first analog of the hypergeometric function for a 
finite place of $\F (x)$''). For $n\in \mathbb Z_+$, $\alpha \in \mathbb 
Z$, denote $D_n=[n][n-1]^q\cdots [1]^{q^{n-1}}$, $[n]=x^{q^n}-x$,
$L_n=[n][n-1]\cdots [1]$ ($n\ge 1$), $D_0=L_0=1$, 
\begin{equation}
(\alpha )_n=\begin{cases}
D_{n+\alpha -1}^{q^{-(\alpha -1)}}, & \text{if $\alpha \ge 1$;}\\
(-1)^{n-\alpha}L_{-\alpha -n}^{-q^n}, & \text{if $\alpha \le 0,n\le -\alpha$;}\\
0, & \text{if $\alpha \le 0,n>-\alpha$}.\end{cases}
\end{equation} 
For $\alpha_i,\beta_i\in \mathbb Z$, such that the series below 
makes sense, Thakur's hypergeometric function is defined as
\begin{equation}
{}_rF_s(\alpha_1,\ldots ,\alpha_r;\beta_1,\ldots 
,\beta_s;z)=\sum\limits_{n=0}^\infty \frac{(\alpha_1)_n\cdots 
(\alpha_r)_n}{(\beta_1)_n\cdots (\beta_s)_nD_n}z^{q^n}
\end{equation}  

Thakur has carried out a thorough investigation of the function 
(1.3) and obtained analogs of many classical properties. In 
particular, he found a differential equation for the function 
(1.3) (see \cite{Th1}). Using the commutation relation $d\tau 
-\tau d=[1]^{1/q}I$, where $I$ is the identity operator, we can 
write the equation from \cite{Th1} in the form
\begin{equation}
\left\{ \prod\limits_{i=1}^r(\Delta -[-\alpha_i])-\left(
\prod\limits_{j=1}^s(\Delta -[-\beta_j])\right) d\right\} 
{}_rF_s=0
\end{equation} 
convenient for us in a sequel. For $r=2,s=1$, this can be written 
in a familiar-looking form: the function $u={}_2F_1(\alpha ,\beta ;
\gamma ;z)$, $\alpha ,\beta ,\gamma \in \mathbb Z$, satisfies the 
equation
\begin{equation}
\tau (1-\tau )d^2u+\left\{ \left( [-1]^q+[-\alpha ]+[-\beta 
]\right) \tau -[-\gamma ]\right\} du-[-\alpha ][-\beta ]u=0.
\end{equation} 
It must be remembered that everywhere in this theory we have to 
deal with nonlinear equations: the operators $\tau$ and $d$ are 
only $\F$-linear, not $K$-, nor $\K$-linear.

A natural next step in developing analysis over $K$
is to try to consider partial differential 
equations with Carlitz derivatives. However, here we encounter a 
serious difficulty noticed in \cite{Khol}: the Carlitz derivatives 
with respect to different variables do not commute. Considering 
the Carlitz rings of ``differential operators'' in the above sense
in \cite{Khol}, the author found a class of partial differential    
operators (acting on an appropriate class of functions of several 
variables), which nevertheless possess reasonable properties. Operators
from this class contain the derivative $d$ in only one distinguished 
variable, and the linear operator $\Delta$ in every other variable. 
Such operators $d$ and $\Delta$ (in different variables) do not 
commute either but satisfy a simple commutation relation.

In this paper we pursue this idea further, introducing a class of
partial differential equations of the above type, for which an 
analog of the Cauchy problem (the initial value problem for a 
differential equation) with respect to the distinguished 
variable is well-posed. In this regard, such equations may be 
seen as function field analogs of classical evolution equations 
of mathematical physics. On the other hand, it is easy to notice an analogy 
between our new equations and general hypergeometric equations of 
classical analysis (see, for example, \cite{BE,GGR1}; note also a 
different non-Archimedean generalization of the hypergeometric 
function \cite{GGR2} where the functions, defined on a 
non-Archimedean field, are complex-valued). In this wider 
framework, it becomes clear how to extend the definition (1.3) of
Thakur's hypergeometric function to the case of the parameters 
from $\K$, so that the resulting function is locally analytic 
with respect to the parameters and coincides with Thakur's one (up 
to a change of variable), as the parameters are of the form 
$[-n]$, $n\in \mathbb N$. Thus, now Thakur's hypergeometric 
function becomes fully an object of the function field arithmetic, with the 
variable, parameters, and values belonging to $\K$.

The structure of the paper is as follows. In Section 2, we prove 
the well-posedness of the Cauchy problem for the above class of 
partial differential equations. We show (Section 3) that these 
equations generate quasi-holonomic modules (in the sense of 
\cite{Khol}). Section 4 is devoted to the general hypergeometric 
function. Analogs of the contiguous relations for the latter are 
given in Section 5. 

\section{Cauchy Problem}

Denote by $\FF_{n+1}$ ($n\ge 1$) the set of all germs of functions 
of the form
\begin{equation}
u=\sum\limits_{i_1=0}^\infty \ldots \sum\limits_{i_n=0}^\infty
\sum\limits_{m=0}^{\min (i_1,\ldots ,i_n)}c_{m,i_1,\ldots 
,i_n}s_1^{q^{i_1}}\ldots s_n^{q^{i_n}}\frac{z^{q^m}}{D_m}
\end{equation}
where $c_{m,i_1,\ldots ,i_n}\in \K$ are such that all the series 
are convergent on some neighbourhoods of the origin. Let $\hat{ 
\mathcal F}_{n+1}$ be the set of polynomials from $\FF_{n+1}$.

Below $d$ will denote the Carlitz derivative in the variable $z$, 
while $\Delta_j$ will mean the difference operator (1.1) in the 
variable $s_j$. In the action of each operator $d,\Delta_j$ on a 
function from $\FF_{n+1}$ (acting in a single variable) other 
variables are treated as scalars. Obviously, linear operators 
$\Delta_j$ commute with multiplications by scalars: 
$\Delta_j\lambda =\lambda \Delta_j$, $\lambda \in \K$, while 
$d\lambda =\lambda^{1/q}d$. We have also the following commutation 
relations:
\begin{equation}
d\tau -\tau d=[1]^{1/q}I,\quad 
d\Delta_j-\Delta_jd=[1]^{1/q}d,\quad \Delta_j\tau -\tau 
\Delta_j=[1]\tau \quad (j=1,\ldots ,n).
\end{equation}
The Carlitz ring $\AG_{n+1}$ (see \cite{Khol}) is generated 
by $\tau ,d,\Delta_j$ ($j=1,\ldots ,n$) and scalars $\lambda \in 
\K$, and the action of any operator from $\AG_{n+1}$ on 
$\FF_{n+1}$ is well-defined.

Let us consider equations of the form
\begin{equation}
\left\{ P(\Delta_1,\ldots ,\Delta_n)+Q(\Delta_1,\ldots 
,\Delta_n)d\right\} u=0
\end{equation}
where $P,Q$ are non-zero polynomials with coefficients from $\K$. 
We look for a solution $u\in \FF_{n+1}$ of the form (2.1) 
satisfying an ``initial condition''
\begin{equation}
\lim\limits_{z\to 0}z^{-1}u(z,s_1,\ldots ,s_n)=u_0(s_1,\ldots ,s_n)
\end{equation}
where $u_0(s_1,\ldots ,s_n)$ is an $\F$-linear holomorphic 
function on a neighbourhood of the origin. The condition (2.4) 
(similar to the initial conditions for ``ordinary'' differential 
equations with Carlitz derivatives \cite{K00}) means actually that 
the coefficients $c_{0,i_1,\ldots ,i_n}$ of the solution (2.1) are 
prescribed for any $i_1,\ldots ,i_n$.

Below we use the notation $[\infty ]=-x$. Then $[n]\to [\infty ]$, 
as $n\to \infty$.

\medskip
\begin{teo}
Suppose that
\begin{equation}
Q([i_1],\ldots ,[i_n])\ne 0\quad \text{for all }i_1,\ldots 
,i_n=0,1,\ldots ,\infty.
\end{equation}
Then the Cauchy problem (2.3)-(2.4) has a unique solution $u\in 
\FF_{n+1}$.
\end{teo}

\medskip
{\it Proof}. It is easy to see that
$$
\Delta_j\left( s^{q^{i_j}}\right) =\begin{cases}
[i_j]s^{q^{i_j}}, & \text{if $i_j\ne 0$};\\
0, & \text{if $i_j=0$}.\end{cases}
$$
The identity $D_m=[m]D_{m-1}^q$ implies the relation
$$
d\left( \frac{z^{q^m}}{D_m}\right) =\begin{cases}
\frac{z^{q^{m-1}}}{D_{m-1}}, &  \text{if $m\ne 0$};\\
0, &  \text{if $m=0$}.\end{cases}
$$
Therefore for a function (2.1) we get
\begin{multline*}
du=\sum\limits_{i_1=1}^\infty \ldots \sum\limits_{i_n=1}^\infty
\sum\limits_{m=1}^{\min (i_1,\ldots ,i_n)}c_{m,i_1,\ldots 
,i_n}^{1/q}s_1^{q^{i_1-1}}\ldots 
s_n^{q^{i_n-1}}\frac{z^{q^{m-1}}}{D_{m-1}}\\
=\sum\limits_{j_1=0}^\infty \ldots \sum\limits_{j_n=0}^\infty
\sum\limits_{\nu =0}^{\min (j_1,\ldots ,j_n)}c_{\nu +1,j_1+1,\ldots 
,j_n+1}^{1/q}s_1^{q^{j_1}}\ldots s_n^{q^{j_n}}\frac{z^{q^\nu }}{D_\nu 
}.
\end{multline*}
Next,
$$
P(\Delta_1,\ldots ,\Delta_n)u
=\sum\limits_{i_1=0}^\infty \ldots \sum\limits_{i_n=0}^\infty
\sum\limits_{m=0}^{\min (i_1,\ldots ,i_n)}c_{m,i_1,\ldots 
,i_n}P([i_1],\ldots ,[i_n])
s_1^{q^{i_1}}\ldots s_n^{q^{i_n}}\frac{z^{q^m}}{D_m}.
$$

Writing a similar formula for $Q(\Delta_1,\ldots ,\Delta_n)u$ and 
substituting all this into (2.3) we find that
\begin{multline*}
\sum\limits_{i_1=0}^\infty \ldots \sum\limits_{i_n=0}^\infty
\sum\limits_{m=0}^{\min (i_1,\ldots ,i_n)}\biggl\{ c_{m,i_1,\ldots 
,i_n}P([i_1],\ldots ,[i_n]) \\
+c_{m+1,i_1+1,\ldots ,i_n+1}^{1/q}Q([i_1],\ldots 
,[i_n])\biggr\} s_1^{q^{i_1}}\ldots 
s_n^{q^{i_n}}\frac{z^{q^m}}{D_m}=0
\end{multline*}
for arbitrary values of the variables. Hence,
we come to the recursion
\begin{equation}
c_{m+1,i_1+1,\ldots ,i_n+1}=-c_{m,i_1,\ldots ,i_n}^q\left\{ 
\frac{P([i_1],\ldots ,[i_n])}{Q([i_1],\ldots ,[i_n])}\right\}^q,
\quad m\le \min (i_1,\ldots ,i_n).
\end{equation}

Since all the elements $c_{0,i_1,\ldots ,i_n}$ ($i_1,\ldots 
,i_n=0,1,2,\ldots$) are given, from (2.6) we find all the 
coefficients of (2.1).

The set $\{ [i],\ i=0,1,\ldots ,\infty \}$ is compact in $K$. 
Therefore the condition (2.5) implies the inequality
\begin{equation}
\left| Q([i_1],\ldots ,[i_n])\right| \ge \mu >0
\end{equation}
for all $i_1,\ldots ,i_n=0,1,2,\ldots ,\infty$. Note also that 
$|[i]|=q^{-1}$ for any $i$, and
\begin{equation}
\left| c_{0,i_1,\ldots ,i_n}\right| \le Cr^{q^{i_1}+\cdots 
+q^{i_n}}
\end{equation}
for some positive constants $C$ and $r$, since the series for the initial condition
$$
u_0(s_1,\ldots ,s_n)=\sum\limits_{i_1=0}^\infty \ldots \sum\limits_{i_n=0}^\infty
c_{0,i_1,\ldots ,i_n}s_1^{q^{i_1}}\ldots s_n^{q^{i_n}}
$$
converges near the origin.

By (2.7) and (2.8),
$$
\left| c_{1,i_1+1,\ldots ,i_n+1}\right| \le C_1^qr^{q^{i_1+1}+\cdots 
+q^{i_n+1}}
$$
(where $C_1>0$ does not depend on $i_1,\ldots ,i_n$), thus
$$
\left| c_{2,i_1+2,\ldots ,i_n+2}\right| \le C_1^{q^2+q}r^{q^{i_1+2}+\cdots 
+q^{i_n+2}},
$$
and, by induction
$$
\left| c_{l,i_1+l,\ldots ,i_n+l}\right| \le C_1^{q^l+q^{l-1}+\cdots +q}
r^{q^{i_1+l}+\cdots +q^{i_n+l}}
$$
for any $l\ge 0$. This means that for any $j_1,\ldots ,j_n\ge l$,
$$
\left| c_{l,j_1,\ldots ,j_n}\right| \le C_2^{q^l}
r^{q^{j_1}+\cdots +q^{j_n}},\quad C_2>0,
$$
which, together with the equality
$$
\left| D_m\right| =q^{-\frac{q^m-1}{q-1}},
$$
implies the convergence of the series in (2.1) near the origin. 
$\qquad \blacksquare$

\medskip
{\it Remark}. It is easy to generalize Theorem 1 to the case of 
systems of equations, where $P$ and $Q$ are matrices whose 
elements are polynomials of $\Delta_1,\ldots ,\Delta_n$. In this 
case, instead of (2.5) we have to require the invertibility of
$Q([i_1],\ldots ,[i_n])$ for all $i_1,\ldots 
,i_n=0,1,\ldots ,\infty$. In an obvious way, this generalization 
covers also the case of a scalar equation of a higher order in 
$d$.

\medskip
Below we will consider in detail some specific examples of the 
equation (2.3). However, before that we consider a holonomic 
property of such equations, valid in the most general situation, 
even without the assumption (2.5).

\section{Quasi-Holonomic Modules}

The class $\FF_{n+1}$ of functions, among which we looked for a 
solution of (2.3), is basic in the function field counterpart 
\cite{Khol} of the theory of holonomic modules \cite{Cout}. Let us 
recall its principal notions.

Every operator $a\in \AG_{n+1}$ can be written in a unique way 
as a finite sum
\begin{equation}
A=\sum a_{l,\mu,i_1,\ldots ,i_n}\tau^ld^\mu \Delta_1^{i_1}\ldots 
\Delta_n^{i_n}, \quad a_{l,\mu,i_1,\ldots ,i_n}\in \K .
\end{equation}
We introduce a filtration in $\AG_{n+1}$ denoting 
by $\Gamma_\nu$, $\nu \in \mathbb Z_+$, the $\K$-vector space of 
operators (3.1) with $\max \{ l+\mu +i_1+\cdots +i_n\} \le \nu$ 
where the maximum is taken over all the non-zero terms contained in (3.1). 
Then $\AG_{n+1}$ becomes a filtered ring.

Suppose $M$ is a left module over $\AG_{n+1}$ with a filtration $\{ \MG_j\}$, 
that is
$$
\MG_0\subset \MG_1\subset \ldots \subset M,\quad 
M=\bigcup\limits_{j\ge 0}\MG_j,
$$
each $\MG_j$ is a finite-dimensional vector space over $\K$,
and $\Gamma_\nu \MG_j\subset \MG_{\nu +j}$ for any $\nu ,j\in 
\mathbb Z_+$. The filtration is called good if the corresponding 
graded module is finitely generated.

For a good filtration there exist a polynomial 
$\chi \in \mathbb Q[t]$ and a number $N\in \mathbb N$, such that
$$
\dim \MG_s=\sum\limits_{i=0}^s\dim (\MG_i/\MG_{i-1})=\chi 
(s)\text{ for }s\ge N
$$
($\dim$ means the dimension over $\K$). The number $d(M)=\deg \chi$, 
called the (Gelfand-Kirillov or Bernstein) {\it 
dimension} of $M$, does not depend on the choice of a good 
filtration. In general, $d(M)\le n+2$. In some cases $d(M)$ can be 
evaluated explicitly. For example, $d(\hat{\mathcal F}_{n+1})=n+1$, while 
$d(\AG_{n+1})=n+2$, if $\AG_{n+1}$ is considered as a left module 
over itself. In contrast to the classical Bernstein theory (see 
\cite{Cout}), in the positive characteristic situation $d(M)$ can 
be arbitrarily small \cite{Khol}. A module $M$ is called {\it 
quasi-holonomic} if $d(M)=n+1$ (in the survey papers \cite{Ksurv,K06}
the term ``holonomic'' is used). As in the classical case, modules 
associated with basic special functions on $K$ possess this 
property \cite{Khol}. In general, the holonomic property is a kind of a 
``quality certificate'' to distinguish objects (equations, functions
etc) with a reasonable behavior.

If $I$ is a non-zero left ideal in $\AG_{n+1}$, then the left 
$\AG_{n+1}$-module $\AG_{n+1}/I$ can be endowed with a good 
filtration, and
\begin{equation}
d\left( \AG_{n+1}/I\right) \le n+1.
\end{equation}

Returning to the equation (2.3) denote by $R$ the operator in the 
left-hand side:
$$
R=P(\Delta_1,\ldots ,\Delta_n)+Q(\Delta_1,\ldots ,\Delta_n)d
$$
where $P,Q$ are non-zero polynomials. Let $I=\AG_{n+1}R$.

\medskip
\begin{teo}
The module $M=\AG_{n+1}/I$ is quasi-holonomic.
\end{teo}

\medskip
{\it Proof}. Due to (3.2), we have to show only that $d(M)\ge n+1$. 
First we prove two lemmas.

\medskip
\begin{lem}
An operator (3.1) is linear if and only if $a_{l,\mu,i_1,\ldots ,i_n}=0$
for $l\ne \mu$.
\end{lem} 

\medskip
{\it Proof}. Let $\sigma \in \K$. Suppose that $A\sigma =\sigma 
A$, that is
$$
\sum \sigma a_{l,\mu,i_1,\ldots ,i_n}\tau^ld^\mu \Delta_1^{i_1}\ldots 
\Delta_n^{i_n}=\sum a_{l,\mu,i_1,\ldots ,i_n}\sigma^{q^{l-\mu }}
\tau^ld^\mu \Delta_1^{i_1}\ldots \Delta_n^{i_n}.
$$
By the uniqueness of the representation (3.1) \cite{Khol}, we find 
that $\sigma^{q^{l-\mu }}=\sigma$, whenever $a_{l,\mu,i_1,\ldots 
,i_n}\ne 0$. Since $\sigma$ is arbitrary, that is possible if and only 
if $l=\mu$. $\qquad \blacksquare$

\medskip
\begin{lem}
The ideal $I$ does not contain non-zero linear operators.
\end{lem} 

\medskip
{\it Proof}. Using the commutation relations (2.2) we can rewrite 
the representation (3.1) of an arbitrary operator $A\in \AG_{n+1}$ 
in the form
\begin{equation}
A=\sum a'_{l,\mu,i_1,\ldots ,i_n}\Delta_1^{i_1}\ldots 
\Delta_n^{i_n}\tau^ld^\mu ,\quad a'_{l,\mu,i_1,\ldots ,i_n}\in \K.
\end{equation}
Just as in (3.1), the coefficients $a'_{l,\mu,i_1,\ldots ,i_n}$ 
are determined in the unique way, and Lemma 1 remains valid for 
the representation (3.3).

Suppose that an operator $A\in \AG_{n+1}$ is such that $AR$ is 
linear. Let us write (3.3) in the form
$$
A=\sum\limits_{l=0}^N\sum\limits_{\mu =0}^N\alpha_{l\mu }\tau^ld^\mu 
$$
where $\alpha_{l\mu }$ are the appropriate elements of the 
commutative $\K$-algebra $\mathcal D$ (without zero divisors) 
generated by the linear operators $\Delta_1,\ldots ,\Delta_n$. We 
have also $R=\gamma +\delta d$, $\gamma ,\delta \in \mathcal D$, 
$\gamma \ne 0$, $\delta \ne 0$. In this new notation,
$$
AR=\left( \sum\limits_{l=0}^N\sum\limits_{\mu =0}^N\alpha_{l\mu 
}\tau^ld^\mu \right) (\gamma +\delta d).
$$

As an element $\gamma \in \mathcal D$ is permuted with powers of 
$\tau$ and $d$, in additional terms (appearing in accordance with 
(2.2)) the powers of $\tau$ and $d$ respectively are the same, 
while the degrees of elements from $\mathcal D$ (as polynomials of 
$\Delta_1,\ldots ,\Delta_n$) decrease. Therefore
$$
AR=\sum\limits_{l=0}^N\sum\limits_{\mu =0}^N\alpha_{l\mu 
}\gamma'\tau^ld^\mu +\sum\limits_{l=0}^N\sum\limits_{\mu =0}^N\alpha_{l\mu 
}\delta'\tau^ld^{\mu +1}
$$
where $\gamma' ,\delta' \in \mathcal D$, $\gamma'\ne 0,\delta'\ne 
0$, whence
$$
AR=\sum\limits_{l=0}^N\sum\limits_{\nu =1}^N\left( \alpha_{l\nu }\gamma'
+\alpha_{l,\nu -1}\delta'\right) \tau^ld^\nu 
+\sum\limits_{l=0}^N\left( 
\alpha_{l0}\gamma'\tau^l+\alpha_{lN}\delta'\tau^ld^{N+1}\right) .
$$

By Lemma 1,
\begin{align}
\alpha_{l0}& =0,\quad l=1,\ldots ,N;\\
\alpha_{lN}& =0,\quad l=0,1,\ldots ,N.
\end{align}
Considering terms with $l<N,\nu =N$, we find that
$$
\alpha_{lN}\gamma' +\alpha_{l,N-1}\delta' =0,
$$
and (3.5) yields $\alpha_{l,N-1}=0$, $0\le l\le N-1$. Repeating the 
reasoning we obtain that $\alpha_{l\nu }=0$ for $l\le \nu$.

On the other hand, for $l\ge 2,\nu =1$ we get
$$
\alpha_{l1}\gamma' +\alpha_{l0}\delta' =0,
$$
and, by (3.4), $\alpha_{l1}=0$. Repeating we come to the 
conclusion that $\alpha_{l\nu }=0$ for $l>\nu$, so that $A=0$. 
$\qquad \blacksquare$

\medskip
{\it Proof of Theorem 2 (continued)}. Let us consider the induced 
filtration in $M$ (see \cite{Cout}). The subspace $\MG_\nu$ is 
generated by images in $M$ of the elements (3.1) with $\max (l+\mu 
+i_1+\cdots +i_n)\le \nu$; those two elements whose difference 
belongs to $I$ are identified. Let us consider elements with 
$l=\mu$.

Elements of the form $\tau^ld^l\Delta_1^{i_1}\cdots 
\Delta_n^{i_n}$ with different collections of parameters 
$(l,i_1,\ldots ,i_n)$ are linearly independent in $\AG_{n+1}$. If 
some linear combination of their images equals zero in $M$, then 
the corresponding linear combination of the elements themselves 
must belong to $I$, which (by Lemma 2) is possible only if it is 
equal to zero. Therefore the images of the above elements are 
linearly independent, so that
\begin{multline*}
\dim \MG_\nu \ge \card \left\{ (l,i_1,\ldots ,i_n)\in \mathbb 
Z_+^{n+1}:\ 2l+i_1+\cdots +i_n\le \nu \right\} \\ \ge
\card \left\{ (l,i_1,\ldots ,i_n)\in \mathbb 
Z_+^{n+1}:\ l+i_1+\cdots +i_n\le [\nu /2]\right\} 
\end{multline*}
(in contrast to the rest of the paper, here $[\cdot ]$ means the 
integer part of a real number).

Evaluating the number of non-negative integral solutions of the 
above inequality as in the proof of Lemma 1 from \cite{Khol} we 
find that
$$
\dim \MG_\nu \ge \dbinom{[\nu /2]+n+1}{n+1} \ge c_1[\nu 
/2]^{n+1}\ge c_2\nu^{n+1}
$$
for large values of $\nu$ ($c_1,c_2$ are positive constants 
independent of $\nu$). Thus, $d(M)\ge n+1$, as desired. $\qquad 
\blacksquare$

\section{Hypergeometric Equations}

Let $n\ge \max (r,s)$, $r,s\in \mathbb N$,
\begin{align*}
P(t_1,\ldots ,t_n) & =\prod\limits_{i=1}^r(t_i-a_i),\\
Q(t_1,\ldots ,t_n) & =\prod\limits_{j=1}^s (t_j-b_j),
\end{align*}
where $a_i,b_j\in \K$, and the elements $b_j$ do not coincide with 
any of the elements $[\nu ]$, $\nu =0,1,\ldots ,\infty$. Then the 
condition (2.5) is satisfied, and the Cauchy problem for the 
equation (2.3) is well-posed.

Let us specify the initial condition (in terms of prescribing the 
values of $c_{0,i_1,\ldots ,i_n}$) as follows:
$$
c_{0,0,\ldots ,0}=1,\quad c_{0,i_1,\ldots ,i_n}=0
$$
for all other values of $i_1,\ldots ,i_n$. Then $c_{m,i_1,\ldots 
,i_n}=0$ for all sets of indices $(m,i_1,\ldots ,i_n)$ except 
those with $m=i_1=\ldots =i_n$. Denote $\sigma_m=c_{m,m,\ldots 
,m}$. By (2.6), we find that
\begin{align*}
\sigma_1 & =\left\{ \frac{(-1)^r\prod\limits_{i=1}^r a_i}{(-1)^s\prod\limits_{j=1}^s 
b_j}\right\}^q,\\
\sigma_2 & =\left\{ \frac{(-1)^r\prod\limits_{i=1}^r a_i}{(-1)^s\prod\limits_{j=1}^s 
b_j}\right\}^{q^2}\left\{ \frac{\prod\limits_{i=1}^r([1]-a_i)}{\prod\limits_{j=1}^s 
([1]-b_j)}\right\}^q,
\end{align*}
and, by induction, after rearranging the factors we get
\begin{equation}
\sigma_m=\frac{ \prod\limits_{i=1}^r 
([0]-a_i)^{q^m}([1]-a_i)^{q^{m-1}}\cdots ([m-1]-a_i)^q}{\prod\limits_{j=1}^s 
([0]-b_j)^{q^m}([1]-b_j)^{q^{m-1}}\cdots ([m-1]-b_j)^q},\quad 
m=1,2,\ldots .
\end{equation}

For $a\in \K$, denote $\langle a\rangle_0=1$,
\begin{equation}
\langle a\rangle_m=([0]-a)^{q^m}([1]-a)^{q^{m-1}}\cdots 
([m-1]-a)^q,\quad m\ge 1
\end{equation}
(of course, $[0]=0$, but we maintain the symbol $[0]$ to have an 
orderly notation). The Pochhammer type symbol $\langle \cdot 
\rangle_m$ satisfies the recurrence
\begin{equation}
\langle a\rangle_{m+1}=([m]-a)^q\langle a\rangle_m^q.
\end{equation}

If $a=[-\alpha ]$, $\alpha \in \mathbb Z$, then
$$
([m]-a)^q=\left( x^{q^m}-x^{q^{-\alpha }}\right)^q=
\left( x^{q^{m+\alpha }}-x\right)^{q^{-(\alpha -1)}}=[m+\alpha ]^{q^{-(\alpha 
-1)}},
$$
so that in this case the recurrence (4.3) coincides with the 
recursive relation \cite{Th1} for the Pochhammer-Thakur symbols 
(1.2). Our normalization $\langle a\rangle_0=1$ is different from 
(1.2) and resembles the classical one.

In the above situation, it follows from (4.1) that the solution 
$u(z;t_1,\ldots ,t_n)$ of the Cauchy problem (2.3)-(2.4) is given 
by the formula
\begin{equation}
u(z;t_1,\ldots ,t_n)={}_rF_s(a_1,\ldots ,a_r;b_1,\ldots 
,b_s;t_1\cdots t_nz)
\end{equation}
where ${}_rF_s$ is the new hypergeometric function
\begin{equation}
{}_rF_s(a_1,\ldots ,a_r;b_1,\ldots 
,b_s;z)=\sum\limits_{m=0}^\infty \frac{ \langle a_1\rangle_m 
\cdots \langle a_r\rangle_m}{ \langle b_1\rangle_m 
\cdots \langle b_s\rangle_m}\frac{z^{q^m}}{D_m}.
\end{equation}
The series (4.5) has a positive radius of convergence since 
$b_1,\ldots ,b_s$ do not coincide with any of the elements $[\nu 
]$, $\nu =0,1,\ldots ,\infty$ (such parameters will be called {\it 
admissible}).

If $a_i=[-\alpha_i]$, $b_j=[-\beta_j]$ ($i=1,\ldots ,r;\ 
j=1,\ldots ,s$), $\alpha_i\in \mathbb Z$, $\beta_j\in \mathbb N$, 
then the function (4.5) coincides with Thakur's hypergeometric 
function (1.3) up to a change of variable $z\Rightarrow \rho z$ 
where $\rho$ depends on all the parameters but does not depend on 
$z$.

Let $h_m$ be the coefficients of the power series (4.4), that is
$$
h_m= \frac{ \langle a_1\rangle_m 
\cdots \langle a_r\rangle_m}{ \langle b_1\rangle_m 
\cdots \langle b_s\rangle_mD_m},\quad m=0,1,2,\ldots .
$$
Since $\dfrac{D_{m+1}}{D_m^q}=[m+1]=\left( 
x^{q^m}-x^{q^{-1}}\right)^q=([m]-[-1])^q$, we find that
\begin{equation}
\frac{h_{m+1}}{h_m^q}=\left\{ \frac{([m]-a_1)\cdots ([m]-a_r)}{([m]-b_1)\cdots 
([m]-b_s)([m]-[-1])}\right\}^q.
\end{equation}

The identity (4.6) means that the ratio $\dfrac{h_{m+1}}{h_m^q}$ 
is the $q$-th power of a rational function of $[m]$, which is a 
clear analog of the basic property of the classical hypergeometric 
function. Note that any rational function of $[m]$ may appear in 
(4.6), except those for which (4.6) does not make sense. Thakur's
hypergeometric function corresponds to the case of rational 
functions with zeroes and poles of the form $[\nu ],\nu \in 
\mathbb Z$.

It follows from (4.4) that $\Delta_iu=\Delta u$ where $\Delta 
=\tau d$ is the difference operator (1.1) in the variable $z$. 
Therefore ${}_rF_s$ satisfies the equation
$$
\left\{ \prod\limits_{i=1}^r\left( \Delta -a_i\right) -
\left( \prod\limits_{j=1}^s\left( \Delta -b_j\right) \right) 
d\right\} {}_rF_s=0.
$$
In particular, for the Gauss-like hypergeometric function 
${}_2F_1$ we have
$$
\{ (\Delta -a)(\Delta -b)-(\Delta -c)d\} {}_2F_1=0.
$$
Substituting $\Delta =\tau d$ and using the commutation relations 
(2.2) we can rewrite this equation in the form
$$
\left\{ \tau (1-\tau )d^2-\left( c-([1]^{1/q}+a+b)\tau \right) 
d-ab\right\} {}_2F_1=0,
$$
which corresponds to (1.5).

If $b\ne [\nu ]$, $\nu =0,1,\ldots ,\infty$, then, in particular, 
$|b+x|\ge \mu >0$, whence
$$
|b-[\nu ]|=\left| (b+x)-x^{q^\nu }\right| =|b+x|\ge \mu
$$
for large values of $\nu$. This means that
$$
|b-[\nu ]|\ge \mu_1>0,\quad \nu =0,1,2,\ldots ,
$$
so that
$$
\left| \langle b\rangle_\nu \right| \ge \mu_1^{q^\nu +q^{\nu 
-1}+\cdots +q}=\mu_1^{\frac{q^{\nu 
+1}-1}{q-1}-1}=\mu_1^{-(q-1)^{-1}-1}\left( 
\mu_1^{\frac{q}{q-1}}\right)^{q^\nu }.
$$
Therefore, if $|z|$ is small enough, then the series (4.5) 
converges uniformly with respect to the parameters $a_i\in \K$ and 
$b_j\in \K \setminus \{ [\nu ],\nu =0,1,\ldots ,\infty \}$ on any 
compact set. Thus, the function ${}_rF_s$ is a locally analytic 
function of its parameters.

\section{Contiguous Relations}

Among many identities for classical hypergeometric functions, a 
special role belongs to relations between contiguous functions, 
that is the hypergeometric functions ${}_2F_1$ whose parameters 
differ by $\pm 1$ \cite{BE}. For Thakur's hypergeometric 
functions, analogs of the contiguous relations were found in 
\cite{Th1,Th2,Th3}. Here we present analogs of the contiguous 
relations for our more general situation. Of 15 possible 
relations, we give, just as Thakur did, only two, which is 
sufficient to demonstrate specific features of the function field 
case; other relations can be obtained in a similar way. Note that 
the specializations of our identities for the case of parameters 
like $[-\alpha ]$, $\alpha \in \mathbb Z$, are slightly different 
from those in \cite{Th1,Th2,Th3}, due to a different normalization 
of our Pochhammer-type symbols.

Our first task is to find an appropriate counterpart of the shift 
by 1 for parameters from $\K$. Denote
\begin{equation}
T_1(a)=(a-[1])^{1/q},\quad a\in \K.
\end{equation}

If $a=[-\alpha ]$, $\alpha \in \mathbb Z$, then
$$
T_1([-\alpha ])=([-\alpha ]-[1])^{1/q}=x^{q^{-\alpha 
-1}}-x=[-\alpha -1],
$$
so that the transformation $T_1$ indeed extends the unit shift of 
integers. The inverse transformation is given by
\begin{equation}
T_{-1}(a)=a^q+[1], \quad a\in \K.
\end{equation}

\medskip
\begin{teo}
The following identities for the Pochhammer-type symbol and the 
Gauss-type hypergeometric function are valid for every $m\in 
\mathbb N$ and any admissible parameters from $\K$:
\begin{gather}
\langle T_1(a)\rangle_m=a^{-q^m}(a-[m])\langle a\rangle_m,\quad 
a\ne 0;\\
\langle a\rangle_{m+1}=-a^{q^{m+1}}\langle T_1(a)\rangle_m^q;\\
\langle T_{-1}(a)\rangle_m=-([1]+a^q)^{q^m}\langle 
a\rangle_{m-1}^q;\\
\langle T_{-1}(a)\rangle_m=-\frac{([1]+a^q)^{q^m}}{([m-1]-a)^q}\langle 
a\rangle_m;
\end{gather}
\begin{equation}
{}_2F_1(T_1(a),b;c;az)-{}_2F_1(a,T_1(b);c;bz)=(a-b){}_2F_1(a,b;c;z);
\end{equation}
\begin{multline}
{}_2F_1(a,b;c;z)-{}_2F_1(a,b;c;z)^q+(c^q-b^q){}_2F_1(a,b;T_1(c);c^{-1}z)^q\\
-(a^q+[1]){}_2F_1(T_{-1}(a),b;c;(a^q+[1])^{-1}z)=0.
\end{multline}
\end{teo}

\medskip
{\it Proof}. Substituting (5.1) into (4.2) and using the fact that 
$[\nu ]+[1]^{1/q}=x^{q^\nu }-x^{1/q}=[\nu +1]^{1/q}$, we get
\begin{equation}
\langle T_1(a)\rangle_m=([1]-a)^{q^{m-1}}([2]-a)^{q^{m-2}}\cdots 
([m]-a),
\end{equation}
which implies (5.3). If we raise both sides of (5.9) to the power 
$q$ and compare the resulting identity with (4.3), we come to 
(5.4). The proofs of (5.5) and (5.6) are similar, based on the 
identity $[\nu ]-[1]=[\nu -1]^q$.

Using (5.3) we find that
$$
{}_2F_1(T_1(a),b;c;az)=\sum\limits_{m=0}^\infty 
(a-[m])\frac{\langle a\rangle_m\langle b\rangle_m}{\langle 
c\rangle_mD_m}z^{q^m},
$$
$$
{}_2F_1(a,T_1(b);c;bz)=\sum\limits_{m=0}^\infty 
(b-[m])\frac{\langle a\rangle_m\langle b\rangle_m}{\langle 
c\rangle_mD_m}z^{q^m},
$$
which implies (5.7). Similarly, if we write down all the terms 
involved in (5.8) and use the identities (4.3), (5.3), (5.4), and (5.6),
after rather lengthy but quite elementary calculations we verify the
required identity (5.8). $\qquad \blacksquare$


\bigskip
\begin{thebibliography}{999}
\bibitem{BE}
H. Bateman and A. Erd\'elyi, {\it Higher Transcendental Functions, 
Vol. 1}, McGraw-Hill, New York, 1953.
\bibitem{Bav}
V. Bavula, The Carlitz algebras, math.RA/0505397.
\bibitem{C35}
L. Carlitz, On certain functions connected with polynomials in a
Galois field, {\it Duke Math. J.} {\bf 1} (1935), 137--168.
\bibitem{Cout}
S. C. Coutinho, {\it A Primer of Algebraic D-modules}, Cambridge 
University Press, 1995. 
\bibitem{GGR1}
I. M. Gelfand, M. I. Graev and V. S. Retakh, General 
hypergeometric systems of equations and series of hypergeometric 
type, {\it Russ. Math. Surv.} {\bf 47}, No. 4 (1992), 1--88.
\bibitem{GGR2}
I. M. Gelfand, M. I. Graev and V. S. Retakh, Hypergeometric 
functions over an arbitrary field, {\it Russ. Math. Surv.} {\bf 59}, No. 5
(2004), 831--905.
\bibitem{G}
D. Goss, {\it Basic Structures of Function Field Arithmetic},
Springer, Berlin, 1996.
\bibitem{K98}
A. N. Kochubei, Harmonic oscillator in characteristic $p$, {\it Lett.
Math. Phys.} {\bf 45} (1998), 11--20.
\bibitem{K99}
A. N. Kochubei, $\F$-linear calculus over function fields, {\it J.
Number Theory} {\bf 76} (1999), 281--300.
\bibitem{K00}
A. N. Kochubei, Differential equations for $\F$-linear functions, {\it J.
Number Theory} {\bf 83} (2000), 137--154.
\bibitem{K03}
A. N. Kochubei, Differential equations for $\F$-linear functions II: Regular 
singularity, {\it Finite Fields Appl.} {\bf 9} (2003), 250--266.
\bibitem{K05}
A. N. Kochubei, Umbral calculus in positive characteristic, {\it 
Adv. Appl. Math.} {\bf 34} (2005), 175--191.
\bibitem{Knlin}
A. N. Kochubei, Strongly nonlinear differential equations with Carlitz 
derivatives over a function field, {\it Ukrainian Math. J.} {\bf 57} (2005), 
794--805.  
\bibitem{Kzeta}
A. N. Kochubei, Polylogarithms and a zeta function for finite 
places of a function field, {\it Contemporary Math.} {\bf 384} (2005), 157--167. 
\bibitem{Khol}
A. N. Kochubei, Quasi-holonomic modules in positive characteristic, 
{\it J. Algebra} {\bf 302} (2006), 826--844.
\bibitem{Ksurv}
A. N. Kochubei, Hypergeometric functions and Carlitz differential equations
over function fields, In: {\it Arithmetic and Geometry around 
Hypergeometric Functions}, Progress in Mathematics, Vol. 260,
Birkh\"auser, Basel, 2007, pp. 163--187.   
\bibitem{K06}
A. N. Kochubei, Umbral calculus and holonomic modules in positive 
characteristic, In: {\it p-Adic Mathematical Physics}, AIP 
Conference Proceedings, Vol. 826 (2006), 254--266.
\bibitem{Sch}
W. Schikhof, {\it Ultrametric Calculus}, Cambridge University
Press, 1984.
\bibitem{Th1}
D. S. Thakur, Hypergeometric functions for function fields, {\it Finite
Fields Appl.} {\bf 1} (1995), 219--231.
\bibitem{Th2}
D. S. Thakur, Hypergeometric functions for function fields II,
{\it J. Ramanujan Math. Soc.} {\bf 15} (2000), 43--52.
\bibitem{Th3}
D. S. Thakur, {\it Function Field Arithmetic}, World Scientific, 
Singapore, 2004.
\end{thebibliography}
\end{document}